\documentclass{amsart}
\usepackage{graphicx}
\newtheorem{theorem}{Theorem}
\theoremstyle{remark}
\newtheorem{remark}[theorem]{Remark}
\numberwithin{equation}{section}
\begin{document}

\title[A short proof of Bing's characterization of $S^3$]{A short proof of Bing's characterization of$S^3$}
\author{Yo'av Rieck}
\address{Department of Mathematical Sciences, 301 SCEN, University of Arkansas,
Fayetteville, AR72701}%
\email{yoav@uark.edu}%


\date{December 21, 2005}

\begin{abstract}
We give a short proof of Bing's characterization of $S^3$: a
compact, connected 3-manifold $M$ is $S^3$ if and only if every
knot in $M$ is isotopic into a ball.
\end{abstract}

\maketitle

Let $M$ be a closed orientable 3-manifold.  We assume familiarity
with the basic notions of irreducible and prime 3-manifolds (see,
{\it e.g.,} \cite{hempel} or \cite{jaco}) and the basic results
about Heegaard splittings of compact 3-manifolds (see, {\it e.g.,
} \cite{scharlemann-review}).  By \em genus \em we always mean
Heegaard genus. A knot $k\subset M$ (that is, a smooth embedding
of the circle into $M$) is called \em irreducible \em if its
exterior $E(k) = M \setminus N(k)$ is an irreducible 3-manifold.
In his own words, Bing's Theorem \cite[Theorem 1]{bing} is:

\begin{theorem}[Bing]
\label{thm}
A compact, connected 3-manifold $M$ is topologically $S^3$ if each %
simple closed curve in $M$ lies in a topological cube in $M$.
\end{theorem}

By ``topological cube" Bing meant what we usually call a ball.
Clearly, any knot in $S^3$ is contained in a ball.  If a knot $k$
in a manifold $M \not\cong S^3$ is contained in a ball (say $B$)
then by considering the boundary of $B$ we see that $k$ is not
irreducible.  Thus, Theorem~\ref{thm} follows from:

\begin{theorem}
\label{thm_too}
Any compact, connected 3-manifold admits an irreducible knot. %
\end{theorem}

%\begin{remark}
In \cite[Theorem 8.1]{jaco-rubinstein} Jaco and Rubinstein gave a
very short proof of Theorem~\ref{thm} for irreducible manifolds,
but their proof relies on the existence of 0-efficient
triangulations.  The purpose of this note is giving a short,
elementary proof of Theorem~\ref{thm}.
%\end{remark}

\noindent{\bf Acknowledgement.}  I would like to thank the referee
for a report that helped make this proof clearer (albeit longer).

\section{The proof}

We prove Theorem~\ref{thm_too}; as remarked above
Theorem~\ref{thm} follows.

\noindent{\bf Case One: $M$ is prime.}  First, when $M$ has genus
at most one, let $k$ be a knot on a Heegaard torus (in $M$) with
$E(k)$ a Seifert fibered space over the disk with 2 exceptional
fibers, which is irreducible.

Second, when $M$ has genus two or more, then $M \not\cong S^2
\times S^1$ and hence is irreducible. Let $M = V_1 \cup_\Sigma
V_2$ be a minimal genus Heegaard splitting of $M$. By Waldhausen
\cite{waldhausen} (see also \cite[Theorem
3.8]{scharlemann-review}) $\Sigma$ is irreducible. Let $k$ be a
core of a 1-handle in $V_1$. Then $\Sigma$ is an irreducible
Heegaard surface for $E(k)$; Haken \cite{haken} (see also
\cite[Theorem 3.4]{scharlemann-review}) showed that every Heegaard
splitting of a reducible manifold is reducible; hence, $E(k)$ is
irreducible.

\begin{remark}
\label{remark}%
In Case One, $\partial E(k)$ is incompressible. For manifolds of
genus one or less this is so by construction of $k$.  For manifold
of genus two or more, if $\partial E(k)$ compressed then (since
$E(k)$ is irreducible) $E(k)$ would be a solid torus; but that
implies $M$ has genus at most one, contradiction.
\end{remark}

\noindent{\bf Case Two: $M$ is composite.}  By Kneser
\cite{kneser} $M$ has a prime decomposition as $M \cong M_1 \#
\cdots \# M_n$ with $M_i$ prime ($i=1,\dots,n$).  Let $k_i \subset
M_i$ be the knot obtained in Case One, let $k = \#_{i=1}^n k_i
\subset M$ be their connected sum, and let $\mathcal{A} \subset
E(k)$ be a collection of annuli that decomposes $k$ into its
summands, that is, the components of $E(k)$ cut open along
$\mathcal{A}$ are homeomorphic to $E(k_i)$ ($i=1,\dots,n)$.

Let $S$ be a sphere in $E(k)$, we will prove that $S$ bounds a
ball. By isotopy of $S$, minimize $S \cap \mathcal{A}$. Assume
that $S \cap \mathcal{A} \neq \emptyset$.  Since $\chi(S)=2$, $S$
cut open along $\mathcal{A}$ has disk components, and let $D$ be
such a disk.  Then $D$ is contained is some component of $E(k)$
cut open along $\mathcal{A}$ (which is homeomorphic to $E(k_i)$,
for some $i$). By Remark~\ref{remark}, $\partial E(k_i)$ is
incompressible and hence $D$ is boundary parallel, contradicting
the minimality assumption.  Hence $S \cap \mathcal{A} =
\emptyset$, and $S$ is contained in a component of $E(k)$ cut open
along $\mathcal{A}$%, which is homeomorphic to $E(k_i)$ (for some
%$i$)
.  By the construction in Case One $S$ bounds a ball. Thus every
sphere in $E(k)$ bounds a ball and $k$ is an irreducible knot,
completing the proof of Theorems~\ref{thm_too} and \ref{thm}.

\providecommand{\bysame}{\leavevmode\hbox
to3em{\hrulefill}\thinspace}
\providecommand{\MR}{\relax\ifhmode\unskip\space\fi MR }
\providecommand{\MRhref}[2]{%
  \href{http://www.ams.org/mathscinet-getitem?mr=#1}{#2}
} \providecommand{\href}[2]{#2}

\end{document}